\documentclass{article}
\usepackage[utf8]{inputenc}

\usepackage{amsmath}
\usepackage{amsthm}
\usepackage{amssymb}
\usepackage{amsfonts}
\usepackage{xcolor}
\usepackage{cancel}
\usepackage{subcaption}
\usepackage[margin=3cm]{geometry}

\usepackage{authblk}

\usepackage{colortbl}

\usepackage{graphicx}
\usepackage{array}
\usepackage{url}

\newcommand{\dom}{\Omega}
\newcommand{\normal}{\boldsymbol{\nu}}

\newcommand{\rot}{\boldsymbol{\nabla} \times }
\newcommand{\Ei}{\mathbf{E}^i}
\newcommand{\E}{\mathbf{E}}

\newcommand{\SDcomment}[1]{}

\title{High Performance Parallel Solvers for the time-harmonic Maxwell Equations}

\date{}
\author[1,2]{Elise Fressart}
\author[2]{Sébastien Dubois}
\author[2]{Loïc Gouarin}
\author[2]{Marc Massot}
\author[1]{Michel Nowak}
\author[2]{Nicole Spillane}

\affil[1]{cortAIx Labs, Thales Research and Technology, 91120 Palaiseau, France}
\affil[2]{CMAP, CNRS, École polytechnique, Institut Polytechnique de Paris, 91120 Palaiseau, France.}

\begin{document}

\maketitle

\paragraph{Abstract}
We consider the numerical solution of large scale time-harmonic Maxwell equations. To this day, this problem remains difficult, in particular because the equations are neither Hermitian nor semi-definite.
Our approach is to compare different strategies for solving this set of equations with preconditioners that are available either in PETSc, MUMPS, or in hypre. Four different preconditioners are considered. The first is the sparse approximate inverse, which is often applied to electromagnetic problems. The second is Restricted Additive Schwarz, a domain decomposition preconditioner. The third is the Hiptmair-Xu preconditioner which is tailored to the positive Maxwell equations, a nearby problem. The final preconditioner is MUMPS's Block Low-Rank method, a compressed block procedure. We also compare the performance of this method to the standard LU factorization technique, which is a direct solver. Performance with respect to the mesh size, the number of CPU cores, the wavelength and the physical size of the domain are considered. This work in progress yields temporary conclusions in favour of the Hiptmair-Xu and the Block Low-Rank preconditioners.

\section{Introduction}

\begin{figure}[h]
     \centering
     \begin{subfigure}[b]{0.45\textwidth}
         \centering
         \includegraphics[width=\textwidth]{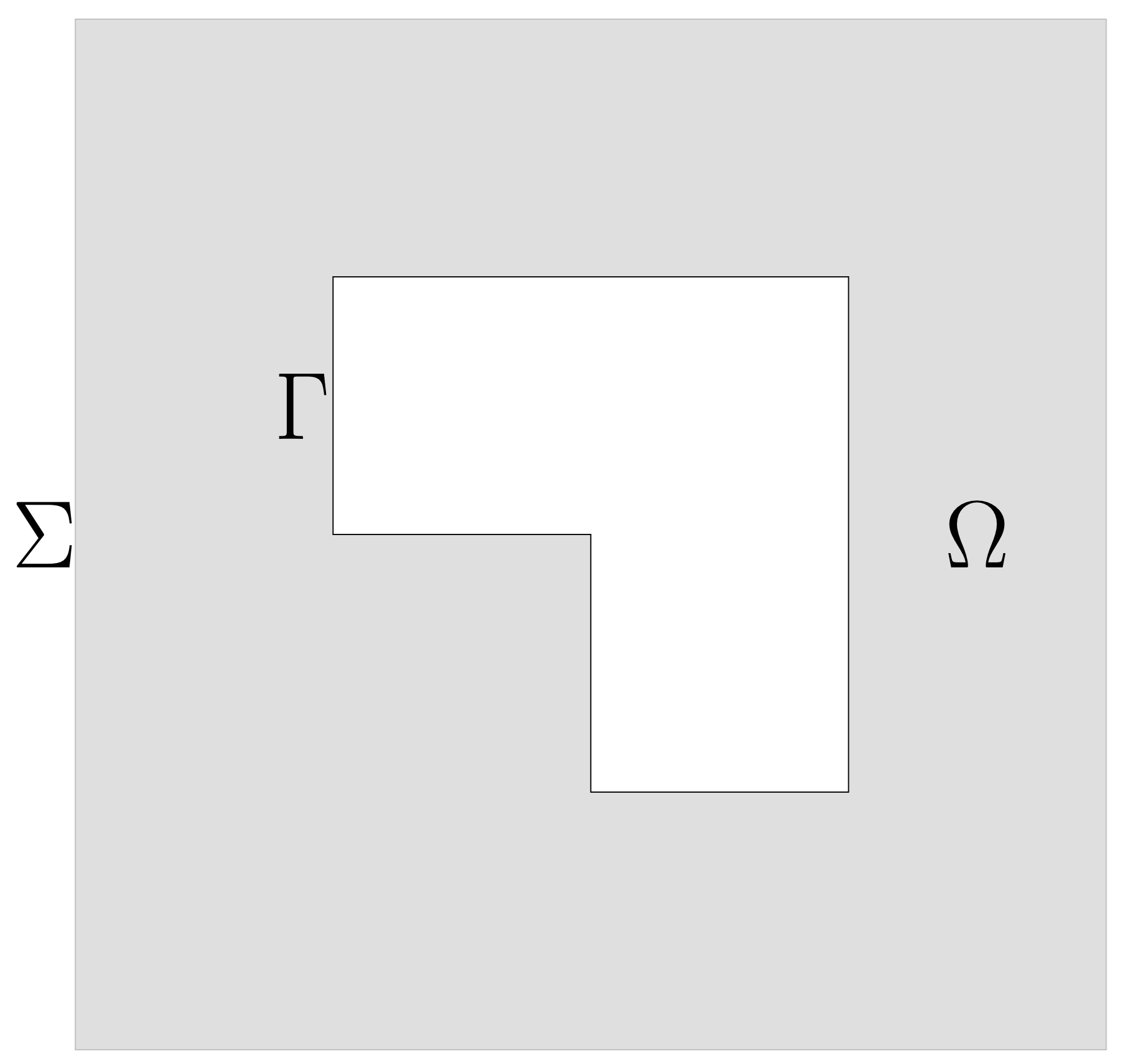}
     \end{subfigure}
     \hfill
     \begin{subfigure}[b]{0.45\textwidth}
         \centering
         \includegraphics[width=\textwidth]{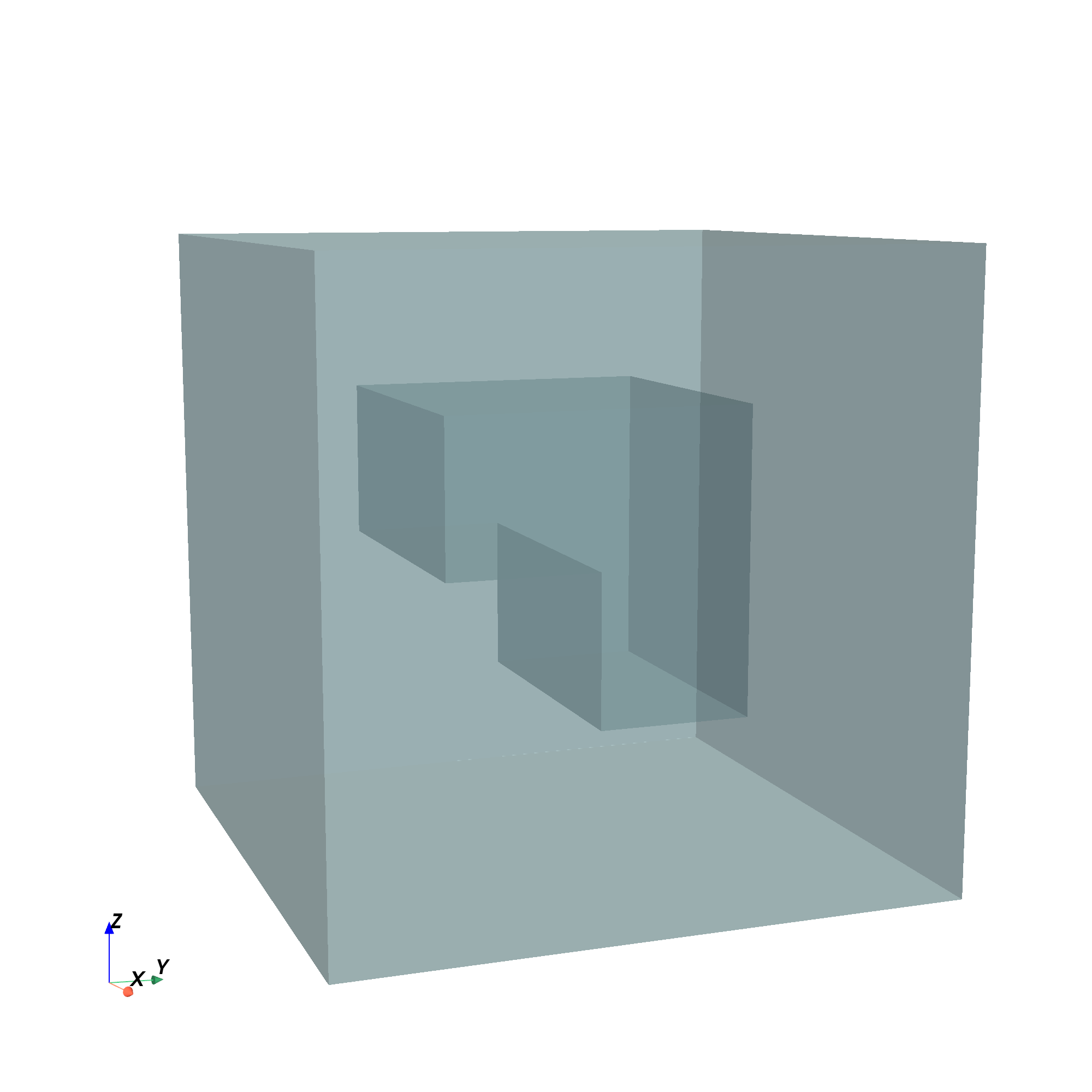}
     \end{subfigure}
        \caption{2D Slice (left) of a 3D (right) possible domain}
        \label{fig:slice_dom_cube_Lshaped}
\end{figure}

We start by introducing the system of equations which we aim to solve numerically: the scattering problem. Let $\dom \subset \mathbb R^2$ or $\mathbb R^3$ be a domain with two boundaries, $\Sigma$ for the exterior domain, and $\Gamma$ for the scatterer (see Figure~\ref{fig:slice_dom_cube_Lshaped}). Denote by $\normal$ the outer unit normal vector on these boundaries. The parameters in the equations are the source term $\mathbf{F}$, the relative permittivity $\epsilon_r$ with the assumption that $\epsilon_r \in \mathbb R$ (lossless medium), the relative permeability $\mu_r$, the impedance $\lambda$ and the incident field $\Ei$ which has wavenumber $k$. With $\E_T = (\normal \times \E) \times \normal$ denoting its tangential component on the boundary $\Sigma$, the scattering problem \eqref{eq:scattering_pb} consists in finding a sufficiently regular electric field $\mathbf{E}$ satisfying
\begin{equation}
\begin{aligned}
    \rot (\mu_r^{-1} \rot \mathbf{E})- k^2 \epsilon_r \mathbf{E} & = \mathbf{F} \quad \text{in }  \dom,\\
    \normal \times \mathbf{E} & = 0 \quad \text{on }  \Gamma, \\
    \mu_r^{-1} (\rot \mathbf{E}) \times \normal - i k \lambda \mathbf{E}_T &= (\rot \Ei) \times \normal - i k \Ei_T \quad \text{on }  \Sigma.
    \label{eq:scattering_pb}
\end{aligned}
\end{equation}

Following \cite[Chapter 4]{zbMATH01868846}, we  next introduce the weak formulation of the problem. We recall the definition of the functional space $\text{H}(\text{curl},\dom)$ and we denote by $X$ the space in which the problem is posed:
\begin{align*}
\text{H}(\text{curl},\dom) & := \left\{ \mathbf{u} \in (\text{L}^2(\dom))^3 ~|~ \rot \mathbf{u} \in (\text{L}^2(\dom))^3 \right\}, \\ 
X & := \left\{ \mathbf{u} \in \text{H}(\text{curl},\dom) ~|~ \normal \times \mathbf{u} = 0 \text{ on }\Gamma \text{ and } \mathbf{u}_T \in (\text{L}^2(\Sigma))^3 \text{ on } \Sigma \right\}.
\end{align*}
 The weak formulation then reads 
\begin{equation}
    \text{Find } \E \in X \text{ such that }  a(\E,\mathbf{u})=b(\mathbf{u}) \quad \forall \mathbf{u} \in X,
    \label{eq:weak_scattering_pb}
\end{equation}
where
\begin{align}
 a(\E, \mathbf{u}) & :=\int_\dom \mu_r^{-1} \rot \E \cdot \rot \bar{\mathbf{u}} - k^2 \int_\dom \epsilon_r \E \cdot \bar{\mathbf{u}} - i k \int_\Sigma \E_T \cdot \bar{\mathbf{u}}_T,  \\
 b(\mathbf{u}) & := \int_\dom \mathbf{F} \cdot \bar{\mathbf{u}} + \int_\Sigma  ((\rot \Ei) \times \normal - i k \Ei_T) \cdot \bar{\mathbf{u}}_T.
 \label{eq:scattering_pb_weak}
\end{align}
Since $\Sigma \neq \emptyset$, this problem has a unique solution \cite[Theorem 4.17]{zbMATH01868846}. The problem is discretized by Nédélec finite elements \cite{zbMATH03653468} of lowest order. The resulting linear system is denoted by
\begin{equation}
    A \E = b,
\end{equation}
in which $\E \in \mathbb C^n$, and the matrix $A \in \mathbb C^{n \times n}$ is both indefinite and non-Hermitian. We have used the notation $\E$ again for the vector of unknowns, this will not cause any confusion since we no longer consider the PDE. 

Common knowledge \cite{zbMATH01953444} is that there are two main families of linear solvers. Direct solvers (such as LU factorization) are robust but, for very large problems, their application becomes impossible. Iterative solvers (such as GMRES) are more natural to parallelize but their convergence may be very slow or not even guaranteed. Accelerating convergence by preconditioning is a necessity for the linear systems that arise from time-harmonic Maxwell. However, finding the correct preconditioner is not an easy task. Some intuition as to why wave propagation problems are so hard to solve with iterative solvers is given in \cite{zbMATH06069657}. Although \cite{zbMATH06069657} is written for the Helmholtz equation, many ideas generalize to Maxwell. 

The purpose of the present work is to explore and compare different available solvers. We present the numerical solution of problems with up to $38$ million degrees of freedom (dofs) using as many as 400 CPU cores on a cluster.  Four preconditioners are considered: the sparse approximate inverse \cite{zbMATH01519797} (Section \ref{sec:sparse_approx_inverse}), one-level overlapping restricted additive Schwarz (RAS) \cite{zbMATH01368737} (Section \ref{sec:ras}), Hiptmair-Xu \cite{zbMATH05485706} (Section \ref{sec:hiptmair-xu}) and LU decomposition with Block Low-Rank (BLR) factorization (in Section \ref{sec:blr}). We also compare these to a direct, full rank, LU factorization. First we setup our problem in Section~\ref{sec:setup}.

But before, we would like to acknowledge that there are other state of the art preconditioners or solvers that we have not yet tested because their implementation is not available to us, or is not adapted to our computational framework. These include multigrid methods \cite{zbMATH02139870,zbMATH07896748}, domain decomposition methods \cite{zbMATH07088069, zbMATH07936343, zbMATH01832364}, sweeping preconditioners \cite{zbMATH07007340, zbMATH06034771}, and $\mathcal{H}-$matrices \cite{10.1007/978-3-031-50769-4_47}.

\section{Computational setup and test cases}
\label{sec:setup}

\paragraph{Hardware and Software}

Our simulations are performed on a cluster with nodes of 2 Intel Xeon Cascade Lake processors (2.10GHz, 20 cores) and with 192 GB RAM by default. This is augmented to 384 GB when needed, in which case we specify it.

 We have implemented our numerical experiments in Python. The linear systems are assembled by the DOLFINx library \cite{Baratta_DOLFINx_the_next_2023}. The lowest order Nédélec elements are named order 1 in DOLFINx. PETSc \cite{petsc-efficient,Dalcin2011} is used to solve these systems. The iterative solver is either left preconditioned GMRES (KSPGMRES in PETSc), or right preconditioned FGMRES \cite{zbMATH00223909}, the flexible GMRES method (KSPFGMRES). Flexible GMRES is adapted to situations where the preconditioner varies throughout the iterations or it is applied in an inexact manner. In both cases, the stopping criterion is a relative residual norm of $10^{-8}$, and both have been applied without restart and with a zero initial guess. The direct solver is the LU decomposition from MUMPS \cite{MUMPS01}.  

MUMPS \cite{MUMPS01} has also been used to obtain the block low-rank factorization in Section~\ref{sec:blr}. For the sparse approximate inverse and the Hiptmair-Xu preconditioner in Sections~\ref{sec:sparse_approx_inverse} and~\ref{sec:hiptmair-xu}, we call the hypre library \cite{hypre-web-page}. In the hypre library, the Hiptmair-Xu preconditioner is referred to as AMS (for  Auxiliary-space Maxwell Solver) \cite{zbMATH05732348}. Finally, the Restricted Additive Schwarz preconditioner in Section \ref{sec:ras} is called directly in PETSc as PCASM.

\paragraph{Split linear system}

We recall that the considered linear system $A \E = b$ arises from discretizing the variational formulation \eqref{eq:weak_scattering_pb} by Nédélec finite elements of the lowest order. The hypre library has the current limitation that it applies only to real valued problems so we rewrite the linear system by splitting the electric field into its real and imaginary parts, $\E_R$ and $\E_I$, respectively. First, we define the matrices associated with the different operators and the right-hand side, with $\{\boldsymbol{\phi}_i ~, 1 \leq i \leq n \}$ denoting the finite element basis:
\begin{align}
C &:= \left( \int_\dom \mu_r^{-1} \rot \boldsymbol{\phi}_i \cdot \rot \boldsymbol{\phi}_j \right)_{1\leq i,j \leq n}, \\
M &:= \left(  k^2 \int_\dom \epsilon_r \boldsymbol{\phi}_i \cdot \boldsymbol{\phi}_j \right)_{1\leq i,j \leq n}, \\
B &:=  \left( k \int_\Sigma \boldsymbol{\phi}_{iT} \cdot \boldsymbol{\phi}_{jT} \right)_{1\leq i,j \leq n}, \\
\mathbf{s}_R &:=  \left( \int_\dom \mathbf{F}_R \cdot \boldsymbol{\phi}_j + \int_\Sigma ((\rot \Ei_R) \times \normal + k \Ei_{I T}) \cdot \boldsymbol{\phi}_{jT} \right)_{1\leq j \leq n}, \\
\mathbf{s}_I &:=  \left( \int_\dom \mathbf{F}_I \cdot \boldsymbol{\phi}_j + \int_\Sigma (\rot \Ei_I)\times \normal - k \Ei_{R T} \cdot \boldsymbol{\phi}_{jT} \right)_{1\leq j \leq n}.
\end{align} 
The split linear system is then
\begin{equation}
\hat A \begin{pmatrix}
\E_R  \\
\E_I \\
\end{pmatrix} =
\begin{pmatrix}
\mathbf{s}_R  \\
-\mathbf{s}_I \\
\end{pmatrix}, \text{ with } \hat A := \begin{pmatrix}
C-M & B  \\
B & -(C-M) \\
\end{pmatrix}.
\label{eq:scattering_problem_split}
\end{equation}

In this notation, the original linear system is 
\begin{equation}
    A \E = b, \text{ with } A := C - M - iB \text{ and } b := s_R + i s_I.
     \label{eq:complex_linear_syst}
\end{equation}

\paragraph{Geometry}
Unless otherwise specified, we solve the problem in the unit cube $[-0.5,0.5]^3$ with an L-shaped scatterer occupying the space $ [-0.25,0.25]^3 \setminus [-0.25,0.25]\times[-0.25,0]\times[-0.25,0]$ (see Fig. \ref{fig:slice_dom_cube_Lshaped}). 

There is one numerical experiment where we dilate the domain by a factor 50 to obtain the $[-25,25]^3$ cube with an L-shaped scatterer occupying $[-12.5,12.5]^3 \backslash [-12.5,12.5]\times[-12.5,0]\times[-12.5,0]$.

\section{Choice of preconditioners and numerical experiments}
\subsection{Sparse approximate inverse}
\label{sec:sparse_approx_inverse}

\paragraph{Presentation of the method} Sparse approximate inverses form a family of often used preconditioners for electromagnetic simulations. A sparse approximate inverse is a matrix $H$ such that the Frobenius norm of $(I - AH)$ is not too large. We tested ParaSails in the hypre library \cite{hypre-web-page} on the split system \eqref{eq:scattering_problem_split}. Following the ideas introduced in \cite{zbMATH01519797}, this algorithm computes a sparsity pattern by sparsifying $A$ (using a parameter \textit{thresh}) and taking the sparsity pattern of the $m$-th power of that matrix. We have set $m=3$ in our experiments. An additional parameter, \textit{filter}, is used in a post-thresholding procedure. An important feature of the sparse approximate inverse, restricted additive Schwarz and Block Low-Rank is that they are purely algebraic meaning that they can be applied without providing anything more to the solver than the matrix $A$ and the right hand side. As the preconditioner remains constant during the iterations, we solve by GMRES.

\paragraph{Analysis of the results from Table~\ref{table:gmres_spai_no_restart_k_10}}
Larger values of \textit{thresh} and \textit{filter} lead to higher number of iterations but smaller solving times to solve the split system \eqref{eq:scattering_problem_split} as can be seen in Table \ref{table:gmres_spai_no_restart_k_10}. For example, if \textit{filter} equals 0.01 and \textit{thresh} equals 0.001, GMRES requires 677 iterations and 19.83s, while it requires 839 iterations and 17.5s for filter equals to 0.05 and threshold equals to 0.01. Small values of \textit{thresh} and \textit{filter} correspond to approximate inverses with a high number of non zero coefficients. With more dofs (not shown here), the number of GMRES iterations reaches the upper limit of 1000. For a fixed set of parameters, the number of iterations increases when the mesh size decreases. Applying the preconditioner to the right does not change significantly the results.

\begin{table}
\begin{subtable}[t]{0.45\textwidth}
\includegraphics[width=0.95\textwidth]{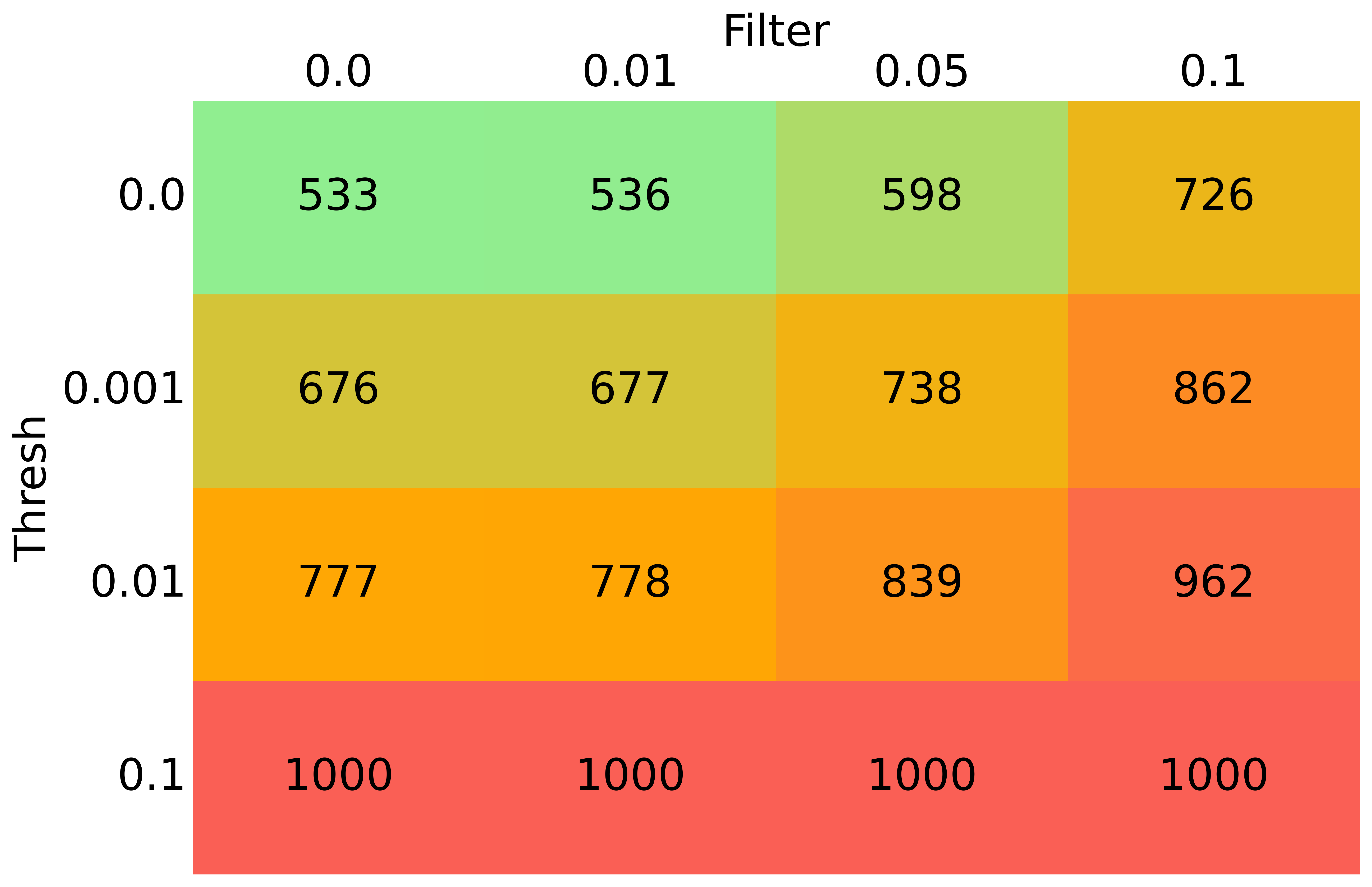}
\end{subtable}
\hfill
\begin{subtable}[t]{0.45\textwidth}
\includegraphics[width=0.95\textwidth]{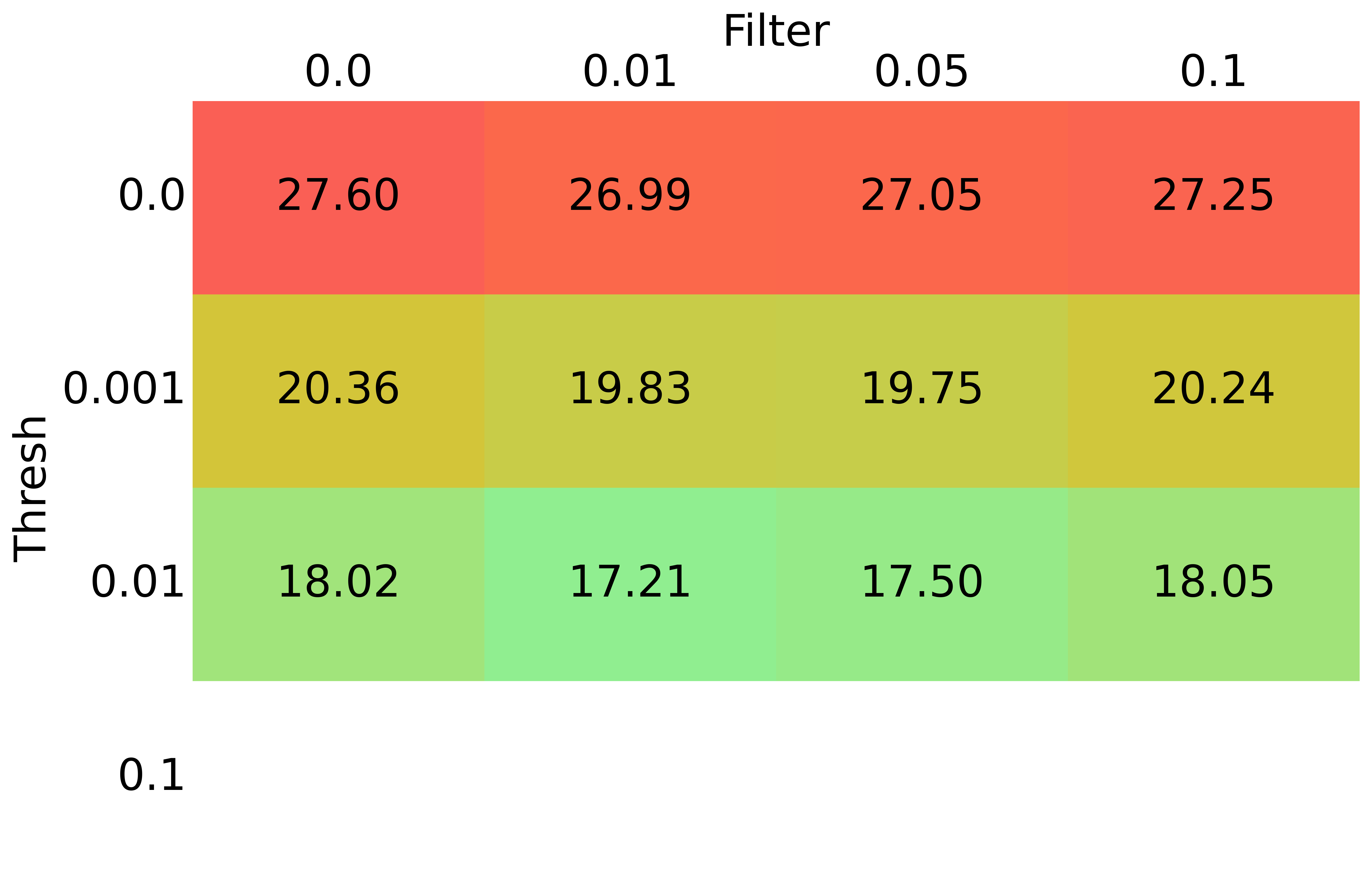}
\end{subtable}
\caption{Number of iterations (left) and solving time (right) for GMRES with a sparse approximate inverse preconditioner (with m=3) for $k=10$ with $10$ points per wavelength (\textit{i.e} $2\times  24,254$ dofs) on 16 CPU cores. The different columns correspond to different filter values. The different rows correspond to different threshold values. When the number of iterations is 1000, the maximum number of iterations is reached while the criterion on the residual is not met. For these cases, the solving time is not shown.}
\label{table:gmres_spai_no_restart_k_10}
\end{table}

\subsection{Restricted Additive Schwarz (RAS)}
\label{sec:ras}

\paragraph{Presentation of the method}
The RAS preconditioner first introduced in \cite{zbMATH01368737} (see also \cite[Sec. 3.9]{zbMATH02113718}) is a domain decomposition preconditioner. The idea of domain decomposition is to partition the simulation domain into subdomains and then approximate the inverse of $A$ by a sum of inverses in each of the subdomains. These are called the subdomain solves. In our experiments, so called exact local solves are used, by performing the LU factorization of the local problems with MUMPS.  It is possible to solve directly the original linear system \eqref{eq:complex_linear_syst} and we do so. 

In detail, the definition of the preconditioner is as follows. The matrix is represented as a graph $G(V,E)$ with vertices $V$ corresponding to the dofs and edges defined by $E=\{(i,j) ~|~ a_{i,j}\neq 0 \}$. The set of vertices is partitioned into $N$ non overlapping subsets denoted $V_i^0$ such that $\cup_{i=1}^N V_i^0=V$. The neighbours of these vertices are then added to form an overlapping subdomain. This is done $\delta$ times to add $\delta$ layers of overlap and the obtained subsets of dofs are denoted by $V_i$. Letting $R_i^0$ and $R_i^\delta$ denote the restriction operators from $V$ to the sets $V_i^0$ and $V_i$, the restricted additive Schwarz preconditioner is defined as
\begin{equation}
    H := \sum_{i=1}^N (R_i^0)^T (R_i^\delta A (R_i^\delta)^T)^{-1} R_i^{\delta}.
\end{equation}

\paragraph{Analysis of the results from Table~\ref{table:gmres_asm_no_restart_k_10}}
We illustrate the difficulty to precondition Maxwell's equations using RAS in Table \ref{table:gmres_asm_no_restart_k_10}. The number of iterations increases with growing number of subdomains (e.g. from 51 with 5 subdomains to 138 with 20 subdomains). The number of iterations decreases when the overlap increases. For the problem with 24,254 dofs and 10 subdomains, it goes from 315 for an overlap of 1 to 34 for an overlap of 5. There is no clear tendency when the mesh is refined. Indeed, with a one-layer overlap refining the mesh doubles the iteration count (e.g. from 315 to 642 on 10 subdomains) whereas with a two-layer overlap, refining the mesh decreases the iteration count (90 versus 58 iterations). 

\begin{table}
\begin{subtable}[t]{0.45\textwidth}
\includegraphics[width=0.95\textwidth]{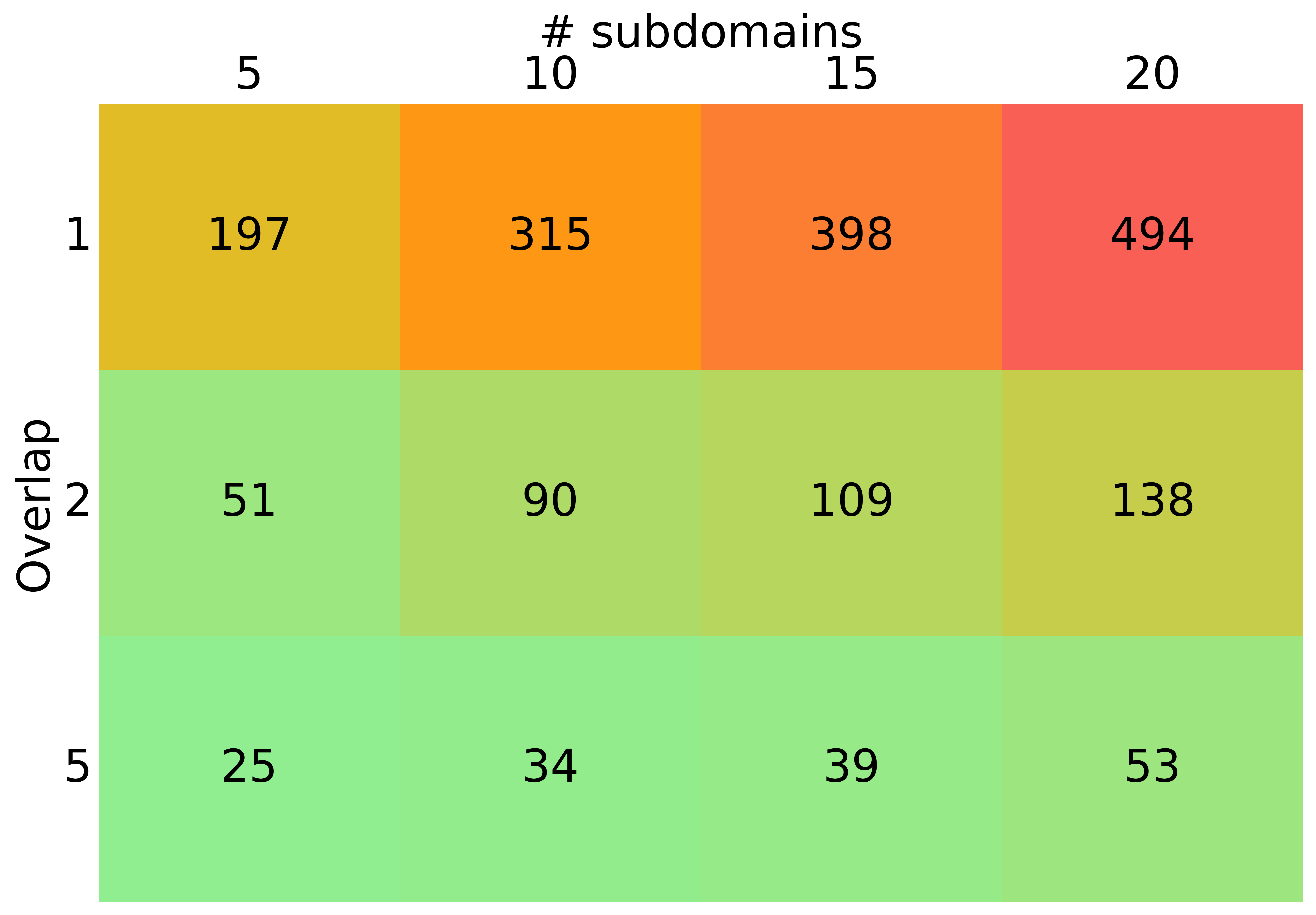}
\end{subtable}
\hfill
\begin{subtable}[t]{0.45\textwidth}
\includegraphics[width=0.95\textwidth]{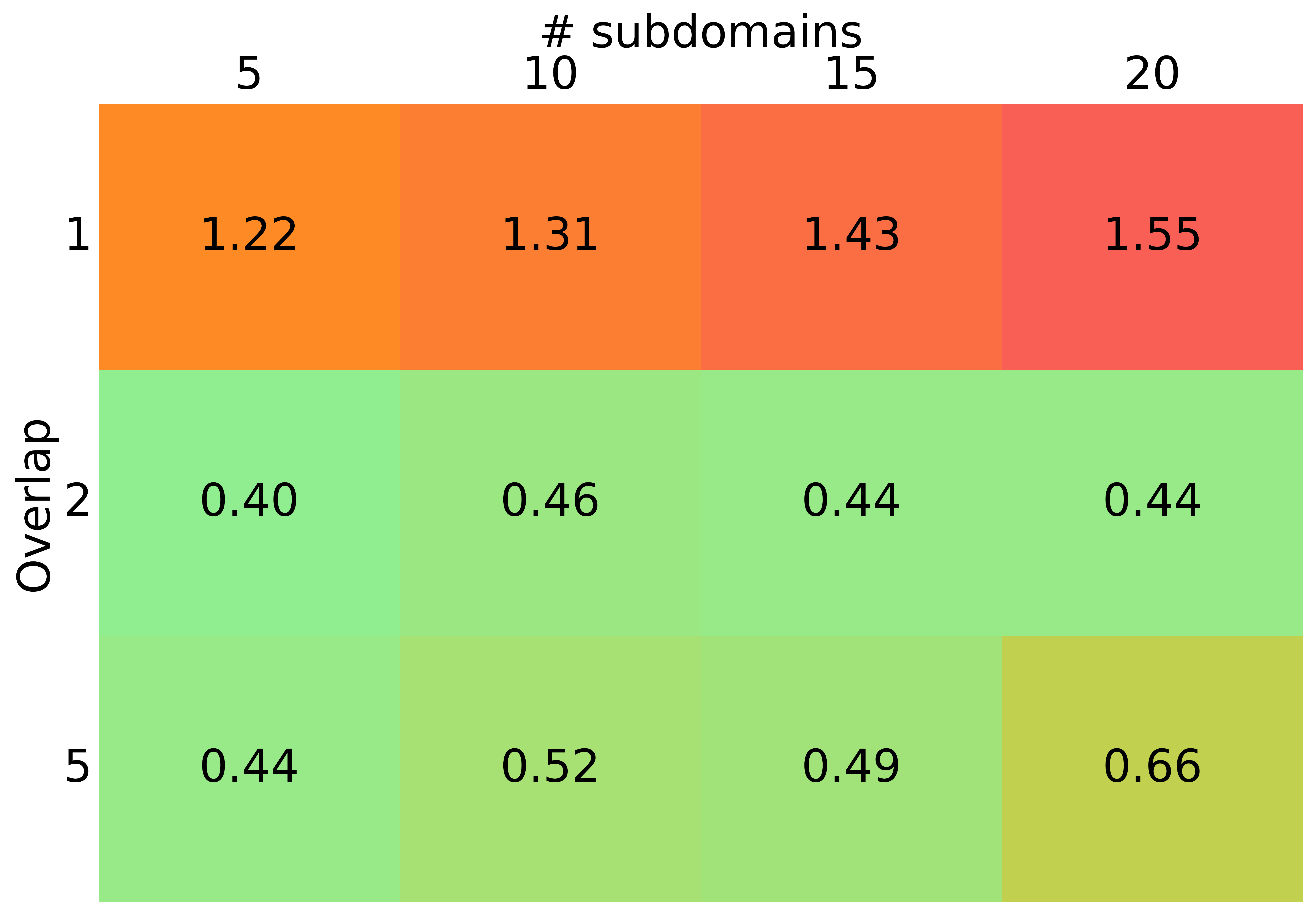}
\end{subtable}

\begin{subtable}[t]{0.45\textwidth}
\includegraphics[width=0.95\textwidth]{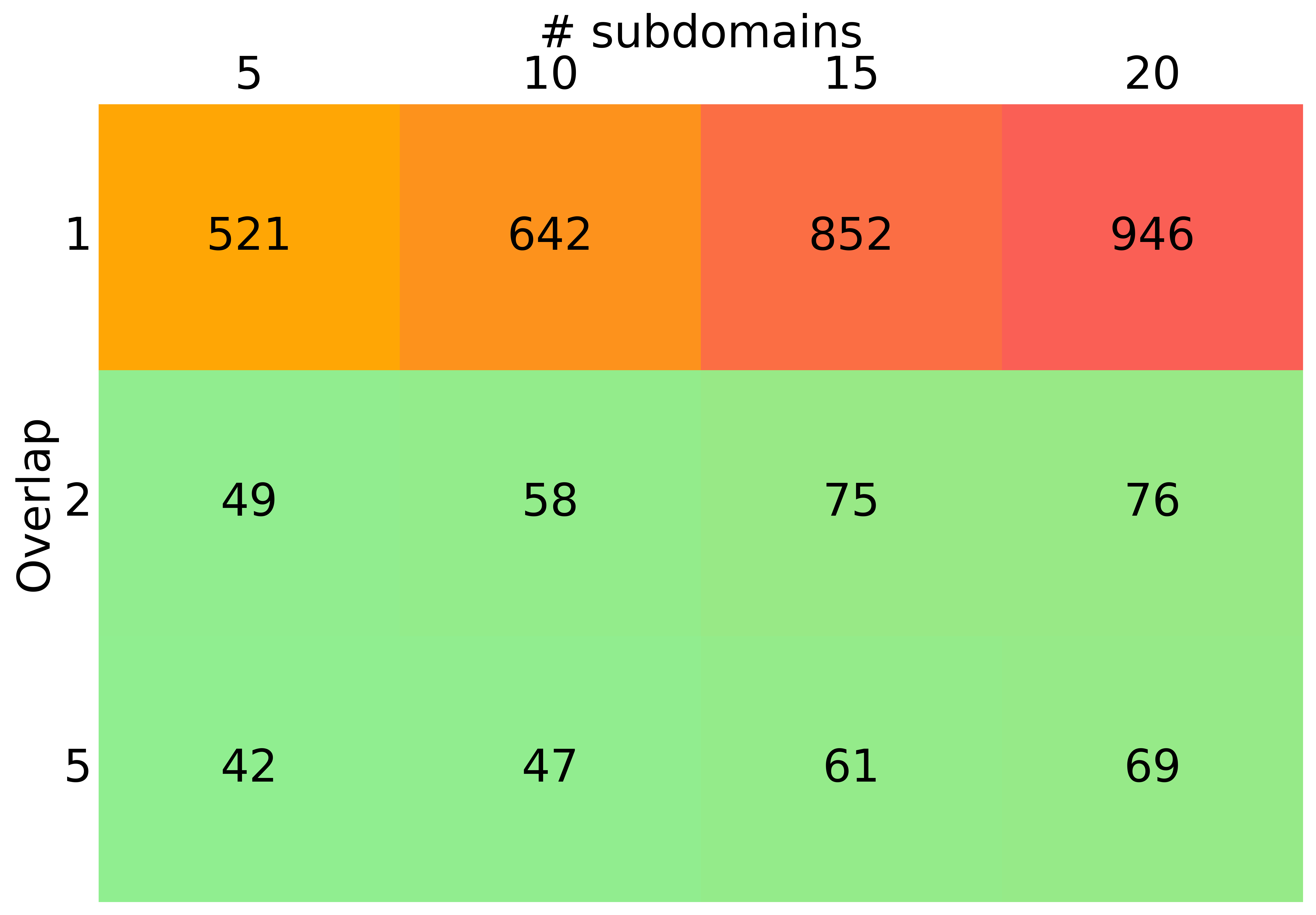}
\end{subtable}
\hfill
\begin{subtable}[t]{0.45\textwidth}
\includegraphics[width=0.95\textwidth]{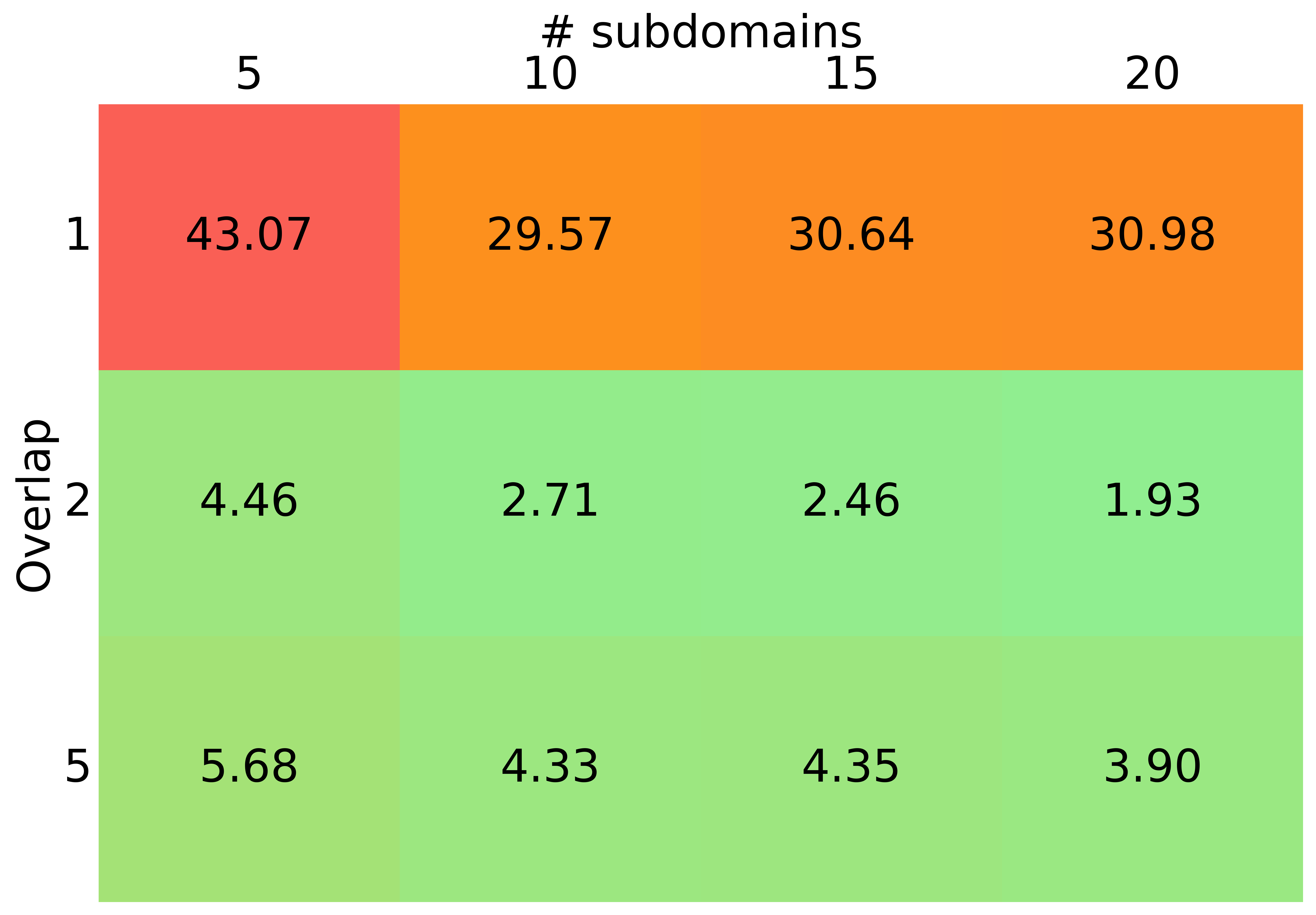}
\end{subtable}

\caption{Number of iterations (left) and solving time (right) for GMRES preconditioned by RAS for $k=10$ with $10$ points per wavelength (\textit{i.e} $24,254$ dofs) (top) and $20$ points per wavelength (\textit{i.e} $172,305$ dofs) (bottom). The different columns correspond to different number of subdomains. The different rows correspond to the number of layers added to each subdomain to build the overlap.}
\label{table:gmres_asm_no_restart_k_10}
\end{table}

\subsection{Iterative solver using Hiptmair-Xu}
\label{sec:hiptmair-xu}

\paragraph{Presentation of the method}
The Hiptmair-Xu (HX) preconditioner \cite{zbMATH05485706} is designed for the positive Maxwell equation with perfect electric conductor boundary conditions, \textit{i.e.}, in its weak formulation
\begin{equation*}
  \text{Find } \E \in \text{H}_0(\text{curl},\dom) \text{ such that } \int_\dom \mu_r^{-1} \rot \E \cdot \rot \bar{\mathbf{u}} + k^2 \int_\dom \epsilon_r \E \cdot \bar{\mathbf{u}} = \int_\dom \mathbf{F} \cdot \bar{\mathbf{u}} \quad \forall \mathbf{u} \in \text{H}_0(\text{curl},\dom),
\end{equation*}
where
\[
\text{H}_0(\text{curl},\dom)  := \left\{ \mathbf{u} \in \text{H}(\text{curl},\dom) ~|~ \normal \times \mathbf{u} = 0 \text{ on }\partial \dom \right\}. 
\]
The Hiptmair-Xu preconditioner is defined by
\begin{equation}
    H_{\text{HX}} := S^{-1} + P_{\text{curl}} (L + k^2 M)^{-1} P_{\text{curl}}^T + (k^2)^{-1} G(-\Delta)^{-1} G^T,
    \label{eq:hiptmairxu_precond}
\end{equation}
where S is a smoother, $L$ is the weighted discrete vector Laplacian matrix on the Lagrangian finite element space, $M$ is the weighted mass matrix on the Lagrangian space, $P_{curl}$ is the matrix of the nodal interpolation operator between the Nédélec space and the Lagrangian space, $G$ is the discrete gradient matrix and $\Delta$ is the weighted discrete scalar Laplacian matrix on the Lagrangian space.

Since the preconditioner is not initially adapted to time-harmonic Maxwell we propose an adaptation inspired by the work of  \cite{grayver2015amssplit,10.1093/gji/ggu119,doi:10.1137/080713112} for time-harmonic Maxwell equations with imaginary relative permittivity, a problem with dissipation, which makes it easier to solve. The authors in \cite{grayver2015amssplit,10.1093/gji/ggu119,doi:10.1137/080713112} propose to precondition the split system (\textit{i.e.} their version of \eqref{eq:scattering_problem_split}) using a block-diagonal matrix where the two diagonal blocks are discretizations of the positive Maxwell equation. The condition number of the preconditioned system is proved to be upper bounded by $\sqrt{2}$ \cite[Lemma 4.1]{doi:10.1137/080713112}. The Hiptmair-Xu preconditioner is then used to invert each block approximately. We adapt this procedure to our problem \eqref{eq:scattering_pb} but we do not recover that the condition number is bounded by $\sqrt{2}$. Similar strategies where the indefinite problem is preconditioned by the positive one are developed, for example, in \cite{Beck1999,zbMATH05240389}. 

We define the following preconditioner $P^{-1}$ for the split problem \eqref{eq:scattering_problem_split} by
\begin{equation}
P := \begin{pmatrix}
C + M + B & 0 \\
0 & C + M + B \\
\end{pmatrix}.
\end{equation}
Instead of factorizing $P$ exactly, each block $(C + M + B )$ is approximately inverted by solving iteratively the corresponding linear systems. For this reason the problem is solved with FGMRES as the outer solver.  For the inner solver (\textit{i.e.}, the one that solves linear systems for $C+M+B$) two methods are compared. The first is to apply HX as a preconditioner: conjugate gradient preconditioned by \eqref{eq:hiptmairxu_precond} is applied with a maximum of 20 iterations and a tolerance of $10^{-2}$ on the relative residual. The second is to apply HX directly as a solver. 

\paragraph{Analysis of the results from Tables~\ref{table:fgmres_ams_no_restart_k_1} to~\ref{table:fgmres_ams_no_restart_k_1_large_domain}}
In the case $k=1$, Table \ref{table:fgmres_ams_no_restart_k_1} compares the time to solve the system and, for the iterative solvers, the number of iterations. The number of points per wavelength varies to study solver performance with respect to the mesh size. Without preconditioner, the GMRES solver converges in fewer than 1000 iterations only for the smallest problem. When Hiptmair-Xu is used as a preconditioner of the inner conjugate gradient, the number of FGMRES iterations does not depend on the mesh size. For these small problems, the direct solver is quicker than the iterative solvers but the difference decreases as the mesh size decreases (from ten times to around three times faster). Simulations running with Hiptmair-Xu as the inner solver are twice as fast as the ones with Hiptmair-Xu as a preconditioner for the finest meshes. However in the former setting, the number of FGMRES iterations is larger and increases from 39 to 52 when the number of dofs increases. 

\begin{table}
\centering
\begin{tabular}{ |c|c|c|c|c|c|c|c|c|}
 \hline
dofs & Direct & \multicolumn{2}{|c|}{No preconditioner} & \multicolumn{3}{|c|}{HX} & \multicolumn{2}{|c|}{HX as solver} \\
 \hline
& time(s) & \#GMRES & time(s) & \#FGMRES & \#CG & time(s) & \#FGMRES & time(s)\\
 \hline
$2\times316$ & $0.01$  & $482$ & $0.61$ &  $28$ & $112+112$ & $0.1$ & $39$ & $0.12$\\
 \hline
 $2\times870$ & $0.02$ &  $>1000^*$  &  $>2.81^*$ & $26$ & $104+104$ & $0.16$ & $40$ & $0.15$\\
   &  &  $(2.4*10^{-2})$  &  & & &  & & \\
 \hline
 $2\times1,440$ & $0.02$ & $>1000^*$ & $>2.83^* $ &  $26$ & $104+104$ & $0.18$ & $42$ & $0.1$\\
 &  & $(1.3*10^{-1})$ & &  & & & & \\
 \hline
 $2\times3,039$ & $0.07$ & $>1000^*$ & $>3.39^* $ &  $26$ & $120+117$ & $0.29$ & $52$ & $0.15$\\
   &  & $(1.7*10^{-1})$ & & & & & & \\
 \hline
 $2\times4,073$ & $0.11$ & $>1000^*$ & $>3.55^* $ &  $26$ & $122+123$ & $0.36$ & $52$ & $0.18$\\
  & & $(1.9*10^{-1})$ & & & & & &\\
 \hline
\end{tabular}
\caption{Solve time and number of iterations for $k=1$ with varying mesh size (10, 20, 30, 40 and 50 points per wavelength) on 8 CPU cores. Superscript * denotes simulations where the maximum number of iterations is reached. For these cases, the final relative norm of the residual is indicated in parentheses.}
\label{table:fgmres_ams_no_restart_k_1}
\end{table}

We study the effect of the wavenumber on the performance of the solvers. As the wavenumber is increased, the problem becomes more indefinite and more difficult to solve. Results for $k=10$ are reported in Table \ref{table:fgmres_ams_no_restart_k_10}. Similar to the case $k=1$, the number of FGMRES iterations does not depend on the mesh size with the Hiptmair-Xu preconditioner. Hiptmair-Xu as a solver is more than twice as fast as the Hiptmair-Xu preconditioner. The number of FGMRES iterations also increases from 190 to 232 iterations with decreasing mesh size. The direct solver is still faster than Hiptmair-Xu used as a preconditioner except for the two largest problems. When it is used as a solver, the iterative solver is faster than the direct solver on meshes with more than 30 points per wavelength. The performance of the Hiptmair-Xu preconditioner deteriorates when $k$ increases (see also Table \ref{table:fgmres_ams_no_restart_k_15}). The number of FGMRES iterations is at least proportional to $k^{0.86}$. For $k=15$ (Table \ref{table:fgmres_ams_no_restart_k_15}) on the two smallest meshes, the direct solver did not run successfully with 300GB of RAM on 32 CPU cores.

\begin{table}
\centering
\begin{tabular}{ |c|c|c|c|c|c|c|}
 \hline
dofs & Direct & \multicolumn{3}{|c|}{HX} & \multicolumn{2}{|c|}{HX as solver} \\
 \hline
& time(s) & \#FGMRES & \#CG & time(s) & \#FGMRES & time(s)\\
 \hline
$2\times24,254$ & $0.67$  & $174$ & $522\times2$ & $4.05$ & $190$ & $1.8$\\
 \hline
 $2\times172,305$ & $7.01$ & $162$ & $646+645$ & $23.23$ & $198$ & $8.52$\\
 \hline
 $2\times556,510$ & $40.38$ & $160$ & $639\times2$ & $75.23$ & $216$ & $30.16$\\
 \hline
 $2\times1,288,926$ & $168.87$ & $154$ & $620+616$ & $168.74$ & $224$ & $73.33$\\
 \hline
 $2\times2,476,448$ & $521.04^{\dagger}$ & $164$ & $717+711$ & $395.28$ & $232$ & $146.08$\\
 \hline
\end{tabular}
\caption{Solve time and number of iterations for $k=10$ with varying mesh size (10, 20, 30, 40 and 50 points per wavelength) on 16 CPU cores. Superscript $\dagger$ indicates that the memory is increased to $300$GB.}
\label{table:fgmres_ams_no_restart_k_10}
\end{table}

\begin{table}
\centering
\begin{tabular}{ |c|c|c|c|c|c|c|}
 \hline
dofs & Direct & \multicolumn{3}{|c|}{HX} & \multicolumn{2}{|c|}{HX as solver} \\
 \hline
& time(s) & \#FGMRES & \#CG & time(s) & \#FGMRES & time(s)\\
 \hline
$2\times76,179$ & $2.89$ & $278$ & $834+834$ & $10.82$ & $294$ & $5.32$\\
 \hline
 $2\times556,510$ & $35.97$ & $264$ & $1052+1053$ & $80.93$ & $308$ & $29.58$\\
 \hline
 $2\times1,823,025$ & $229.97^\dagger$ &  $260$ & $1038+1038$ & $267.05$ & $318$ & $98.73$\\
 \hline
 $2\times4,261,059$ & * &  $254$ & $1016+1016$ & $610.6$ & $334$ & $274.24$\\
 \hline
 $2\times8,239,279$ & * &  $254$ & $1020+1021$ & $1255.95$ & $338$ & $584.66$\\
 \hline
\end{tabular}
\caption{Solve time and number of iterations for $k=15$ with varying mesh size (10, 20, 30, 40 and 50 points per wavelength) on 32 CPU cores. The symbol * indicates that the direct solver failed due to memory issues. Superscript $\dagger$ indicates that the memory is increased to $300$GB.}
\label{table:fgmres_ams_no_restart_k_15}
\end{table}

We also perform the study on the domain enlarged by a factor 50 (see Table \ref{table:fgmres_ams_no_restart_k_1_large_domain}). The preconditioner does not perform well on the enlarged domain even with a small wavenumber. The number of iteration is multiplied by 28 compared the small domain. On 64 CPU cores, the direct solver is twice as fast as the iterative solver with inner preconditioned conjugate gradient. However, the iterative solver succeeds in approximating the solution with 38 million dofs contrary to the direct solver. As expected, the iterative solvers scale better than the direct one. We can compare the scalability of the different solvers. Let $n_{\text{cpu}}$ denote the number of CPU cores and $t$ the solve time. For the direct solver, the solve time is proportional to $n_{\text{cpu}}^{-0.48}$. For the iterative solvers, the scaling is better with $t \propto n_{\text{cpu}}^{-0.83}$ and $t \propto n_{\text{cpu}}^{-0.84}$ for Hiptmair-Xu used as a preconditioner for CG and as the inner solver respectively.

\begin{table}
\centering
\begin{tabular}{ |c|c|c|c|c|c|c|c|}
 \hline
dofs & CPU cores & Direct & \multicolumn{3}{|c|}{HX} & \multicolumn{2}{|c|}{HX as solver} \\
 \hline
& & time(s)& \#FGMRES & \#CG & time(s) & \#FGMRES & time(s)\\
 \hline
$2\times 2,474,860$ & $8$ & $684.94^{\dagger}$ &  $784$ & $2352+2352$ & $2722.41$ & $806$ & $1278.69$\\
 \cline{2-8}
& $16$ & $525.09^{\dagger}$ & $784$ & $2352+2352$ & $1453.66$ & $810$ & $686.29$\\
 \cline{2-8}
& $32$ & $390.0^{\dagger}$ &  $784$ & $2352+2352$ & $986.29$ & $812$ & $473.36$\\
 \cline{2-8}
& $64$ & $252.7$ &  $784$ & $2352+2352$ & $480.8$ & $814$ & $225.06$\\
 \hline
 \hline
$2\times 19,271,610$ & $200$ & * &  $746$ & $2238+2238$ & $1568.7$ & $851$ & $935.3$\\
 \cline{2-8}
& $400$ & * &  $746$ & $2238+2238$ & $821.87$ & $856$ & $485.99$\\
 \hline
\end{tabular}
\caption{Solve time and number of iterations for $k=1$ with 10 and 20 points per wavelength in the enlarged domain with varying number of CPU cores. Superscript $\dagger$ indicates that the memory is increased to $300$GB.}
\label{table:fgmres_ams_no_restart_k_1_large_domain}
\end{table}

\subsection{Block Low-Rank factorization as preconditioner}
\label{sec:blr}

\paragraph{Presentation of the method} 
Although they are tremendously efficient on a very large range of problems, direct solvers face significant limitations associated with their computational complexity ($O(n^2)$) and their substantial memory footprint ($O(n^{4/3})$) \cite{mary:tel-01929478}. Indeed, these characteristics can quickly render them impracticable for the problem sizes encountered in realistic electromagnetic simulations. 

To mitigate these challenges while retaining some of the desirable characteristics of direct methods, Block Low-Rank (BLR) approximations offer a promising approach.
BLR techniques belong to the broader family of Hierarchical Matrix ($\mathcal{H}$-matrix) methods \cite{zbMATH01967713}. The core principle of BLR is to exploit the low numerical rank of the off-diagonal matrix blocks by compressing them. Indeed, many of the dense blocks that appear in factors can be well approximated by low-rank matrices, particularly those corresponding to interactions between degrees of freedom that are `far apart' from each other in the underlying graph. By identifying and compressing these blocks, both memory requirements and computational costs for the factorization and the solution phases can be significantly reduced: the FLOP complexity decreases from $O(n^2)$ in full rank to $O(n^{1.67})$ in BLR, while the memory complexity decreases from $O(n^{4/3})$ to $O(n \log(n))$ \cite{doi:10.1137/16M1077192}. The BLR functionality in MUMPS is primarily controlled by a user-defined parameter $\varepsilon$. This parameter determines the threshold for the low-rank compression applied to eligible frontal blocks during factorization. When a dense block is identified as compressible within the LU factors, a low-rank approximation is computed. A smaller $\varepsilon$ value results in a more accurate low-rank approximation, leading to less compression \cite{doi:10.1137/120903476}. 
We use the resulting compressed matrix as a preconditoner for FGMRES. Our goal is to evaluate the effectiveness of this approach and assess whether BLR can serve as a viable preconditioner in terms of convergence, overall runtime and memory consumption.

\paragraph{Analysis of results from Table~\ref{table:blr_as_direct_residual}}

\begin{table}
\centering
\begin{tabular}{ |c|c|c|c|}
 \hline
$\varepsilon$ & time(s) & \#FGMRES & Compression \\
 \hline
$5.0 \times  10^{-3}$   & 342 & 56  & 2.25  \\ \hline
$10^{-3}$               & 321 & 10  & 2.05  \\ \hline
$10^{-5}$               & 340 & 3   & 1.53  \\ \hline
$10^{-9}$               & 476 & 2   & 1.03  \\ \hline \hline
FR                      & 411$^\dagger$ & -   & -     \\ \hline
\end{tabular}
\caption{Solve time and number of iterations of the FGMRES with varying $\varepsilon$ values on the $2\times 2,474,860$ dofs matrix, using $k=1$ and 10 points per wavelength, with 32 CPU cores. The compression is the ratio of factor sizes between the corresponding compressed factorisation and the full-rank one. The target relative residual is set to $10^{-8}$. Superscript $\dagger$ indicates that the memory is increased to $350$GB.}
\label{table:blr_as_direct_residual}
\end{table}

Table~\ref{table:blr_as_direct_residual} illustrates the trade-off between compression, computational time and memory consumption when varying the BLR compression threshold $\varepsilon$ value. As expected, more compression leads to a significant reduction of the factor size (i.e., $\varepsilon = 10^{-3}$, with a compression ratio of 2.05). Despite the relatively coarse approximation, using this BLR factorization as a preconditioner remains effective, since GMRES converges in only 10 iterations with $\varepsilon = 10^{-3}$. However, the preconditioning quality deteriorates for looser threshold $\varepsilon \geq 5.0 \times 10^{-3}$, as reflected by an increase in the number of GMRES required iterations for convergence. This increase results in a higher walltime, making very aggressive BLR compressions less efficient.

In conclusion, the results show that BLR factorizations not only compresses efficiently the matrices arising from this class of electromagnetics problems,  but also provides effective preconditioners. 

\subsection{Comparison of HX and BLR in Table~\ref{table:blr_hx_large_dom}}

In Table \ref{table:blr_hx_large_dom}, we group results from Tables \ref{table:blr_as_direct_residual} and \ref{table:fgmres_ams_no_restart_k_1_large_domain} to compare the Hiptmair-Xu approach and the Block Low-Rank on the same test case. The BLR preconditioner is faster and requires fewer iterations for this problem. Since RAS did not converge on this test case after $1000$ iterations, results for RAS have not been included in the table. 

\begin{table}
\centering
\begin{tabular}{ |c|c|c|c|}
 \hline
 & time(s) & \#FGMRES \\
 \hline
BLR$(10^{-3})$               & 321 & 10 \\ \hline
BLR$(10^{-5})$               & 340 & 3  \\ \hline \hline
FR                      & 411 $^\dagger$ & -     \\ \hline \hline
HX & 986.29 & 784 \\ \hline 
HX as solver & 473.36 & 812 \\ \hline 
\end{tabular}
\caption{Solve time and number of iterations of the FGMRES, using $k=1$ and 10 points per wavelength (\textit{i.e.} $2\times 2,474,860$), with 32 CPU cores. Superscript $\dagger$ indicates that the memory is increased to $350$GB.}
\label{table:blr_hx_large_dom}
\end{table}

\section{Conclusion}

To conclude, we have compared five high performance solvers for the time-harmonic Maxwell equations with particular emphasis on FGMRES preconditioned either by Hiptmair-Xu, or by a Block Low-Rank preconditioner. Hiptmair-Xu exhibits a convergence independent of the mesh size for the problem considered but its performance deteriorates with increasing wavenumber $k$, and on physically larger subdomains. With BLR, there is a compromise to be found between fast convergence and memory consumption. In future work we plan on pursuing this comparison further, on even larger problems and to include more solvers into the analysis.

\section*{Acknowledgements}

The present work is funded by a CIFRE Grant between Thales and CMAP at \'Ecole polytechnique (E. Fressart) and has been supported by the HPC@Maths Initiative (PI. L. Gouarin and M. Massot), which we would like to acknowledge. The support of the Exa-MA project (Methods and Algorithms for Exascale, part of the PEPR NumPEx, Grant/Award Number: ANR-22-EXNU-0002 - https://numpex.org/) is also gratefully acknowledged. 
This work has been granted access to the HPC resources of IDCS support unit from \'Ecole polytechnique.

\bibliographystyle{unsrt}
\bibliography{bibliography}
\end{document}